\documentclass[10pt]{article}
\usepackage{amsmath}
\usepackage{amssymb}
\usepackage{amsfonts}

\def\noi{{\noindent}}
\def\cq{ \hfill $\blacksquare$ }
\def\llbracket{[\hspace{-.10em} [ }
\def\rrbracket{ ] \hspace{-.10em}]}

\def\ii{{\cal I}}

\def\t{{\cal T}}

\def\r{{\cal R}}
\def\z{{\cal Z}}
\def\n{{\cal N}}

\def\ov{\overline}

\def\la{\longrightarrow}
\def\da{\downarrow}

\def\noi{\noindent}

\def\build#1_#2^#3{\mathrel{\mathop{\kern 0pt#1}\limits_{#2}^{#3}}}

\newtheorem{theorem}{Theorem}[section]
\newtheorem{lemma}[theorem]{Lemma}
\newtheorem{proposition}[theorem]{Proposition}
\newtheorem{corollary}[theorem]{Corollary}
\newtheorem{definition}{Definition}[section]

\newcommand{\R}{\mathbb{R}}
\newcommand{\Z}{\mathbb{Z}}

\newcommand{\N}{\mathbb{N}}

\newcommand{\T}{\mathbb{T}}

\newcommand{\un}{\boldsymbol{1}}

\begin{document}

\title{ \bf RANDOM REAL TREES}
\author{ by \\
\\
Jean-Fran\c cois {\sc Le Gall} \\
{\small D.M.A., Ecole normale sup\'erieure, 45 rue d'Ulm, 75005 Paris, France} \\
{\small\tt legall@dma.ens.fr}}
\vspace{4mm}
\date{}
\maketitle

\begin{abstract}
We survey recent developments about random real trees, whose prototype is
the Continuum Random Tree (CRT) introduced by Aldous in 1991. We briefly explain the
formalism of real trees, which yields a neat presentation of the theory 
and in particular of the relations between discrete Galton-Watson trees
and continuous random trees.
We then discuss the particular class of self-similar random real trees 
called stable trees, which generalize the CRT. We review several important
results concerning stable trees, including their branching property, which is analogous
to the well-known property of Galton-Watson trees, and the calculation 
of their fractal dimension. We then consider
spatial trees, which combine the genealogical structure of a real tree with spatial displacements,
and we explain their connections with superprocesses. In the last section, we
deal with a particular conditioning problem for spatial trees, which is 
closely related to asymptotics for random planar quadrangulations.

\begin{center}
{\bf R\'esum\'e}
\end{center}

Nous discutons certains d\'eveloppements r\'ecents de la th\'eorie des arbres r\'eels
al\'eatoires, dont le prototype est le CRT introduit par Aldous en 1991. Nous introduisons
le formalisme d'arbre r\'eel, qui fournit une pr\'esentation \'el\'egante de 
la th\'eorie, et en particulier des relations entre les arbres de Galton-Watson discrets
et les arbres continus al\'eatoires. Nous discutons ensuite la classe des arbres 
auto-similaires appel\'es arbres stables, qui g\'en\'eralisent le CRT. Nous 
pr\'esentons plusieurs r\'esultats importants au sujet des arbres stables,
notamment leur propri\'et\'e de branchement, analogue continu d'une propri\'et\'e bien
connue pour les arbres de Galton-Watson, et le calcul de leurs dimensions fractales.
Nous consid\'erons ensuite les arbres spatiaux, qui combinent la structure g\'en\'ealogi\-que d'un
arbre r\'eel avec des d\'eplacements dans l'espace, et nous expliquons leurs liens
avec les superprocessus. Dans la derni\` ere partie, nous traitons un conditionnement
particulier des arbres spatiaux, qui est \'etroitement li\'e \` a certains 
r\'esultats asymptotiques pour les quadrangulations planes al\'eatoires.
\end{abstract}

\section*{Introduction}

The purpose of this paper is to give an overview of recent work about
continuous genealogical structures and their applications. The interest for these
continuous branching structures first arose from their connections with the
measure-valued branching processes called superprocesses, which have been studied
extensively since the end of the eighties. Since superprocesses are obtained as
weak limits of branching particle systems, it is not surprising that their
evolution should be coded by a kind of continuous genealogical structure, and
Perkins used non-standard analysis to give a precise definition 
of this structure (see in particular \cite{Per1}).
A little later, the Brownian snake construction of
(finite variance) superprocesses \cite{LG0},\cite{LG1} provided another way of describing the
underlying
genealogy. This construction made it clear that the genealogy  of superprocesses, or equivalently 
of Feller's branching diffusion process (which corresponds to the
total mass of a superprocess),
could be coded by the structure of excursions of linear Brownian motion above
positive levels. The Brownian snake approach had significant applications to
sample path properties of superprocesses \cite{LGP} or to their connections
with partial differential equations \cite{LG3}, \cite{LG4}. 

In a series of papers
\cite{Al1}, \cite{Al2}, \cite{Al3} at the beginning of the nineties, Aldous 
developed the theory of the Continuum Random Tree, and showed that this object 
is the limit as $n\to\infty$, in a suitable sense, of rescaled critical Galton-Watson trees 
conditioned to have $n$ vertices (see Theorem \ref{Aldous} below). Although the CRT was first defined as a 
particular random subset of the space $\ell^1$, it was identified in \cite{Al3}
as the tree coded by a normalized Brownian excursion, in a way very similar to
the Brownian snake approach to superprocesses (note however that the CRT is related
to a Brownian excursion normalized to have duration $1$, whereas in the Brownian snake
approach it is more natural to deal with unnormalized excursions).

In a subsequent paper, Aldous \cite{Al4} suggested the definition of the
so-called integrated super-Brownian excursion (ISE), which combines the
genealogical structure of the CRT with spatial Brownian displacements. A simple way of
looking at
ISE is to view it as the uniform measure on the range of a Brownian snake
driven by a normalized Brownian excursion (see Section IV.6 in \cite{LG4}, and also
Definition 6.1 below). ISE turned 
out to appear in asymptotics for several models of statistical mechanics: See
in particular \cite{Sl1} and \cite{Sl2}. In certain cases however, the continuous
branching structure of ISE (or equivalently of the CRT) is not appropriate to describe
the asymptotics of the model: For instance, the results of \cite{Sl3} for oriented percolation 
involve the canonical measure of super-Brownian motion, whose genealogical struture
is described by a Brownian excursion with height greater than $1$, rather than
a normalized excursion. This suggests that rather than concentrating on the
normalized excursion it is worthwile to deal with various types of Brownian excursions
which correspond to different conditionings of the fundamental object which is the
It\^o measure. 

Both the Brownian snake approach to superprocesses and Aldous' representation of the CRT
correspond to the fact that the genealogical structure of Feller's branching diffusion process 
can be coded by positive Brownian excursions, in a sense that can be made very precise
via the considerations developed in Section 2 below. It was then natural to ask for
a similar coding of the genealogy of more general continuous-state branching processes. Recall that
continuous-state branching processes are Markov processes with values in $\R_+$, which are the possible
limits of rescaled Galton-Watson branching processes (Lamperti \cite{Lam1}).
Such processes are characterized by a branching mechanism function $\psi$, with
$\psi(u)=c\,u^2$ in the case of Feller's diffusion.
The problem of coding the genealogy of general (critical or subcritical)
continuous-state branching processes was treated in two papers of Le Gall and Le
Jan \cite{LGLJ1},
\cite{LGLJ2} (see also the monograph \cite{DuLG0}). The role of the Brownian excursion 
in the case of Feller's diffusion is
played by the so-called height process, which is a (non-Markovian) function of the
L\'evy process with no negative jumps and Laplace exponent $\psi$. The construction of the height process and
its relations with the genealogy of continuous-state branching processes
made it possible to
investigate the asymptotics of critical Galton-Watson trees
when the offspring distribution has infinite variance
(see Chapter 2 of \cite{DuLG0}, and \cite{Duq}). See
Section 4 below for a brief presentation in the stable case where $\psi(u)=u^\alpha$
for some $1<\alpha\leq 2$. 

In the present work, we give a survey of the preceding results, and
of some recent applications, using the language and the formalism of real trees
(cf Section 1).
Although real trees have been studied for a long time for algebraic or geometric purposes
(see e.g. \cite{DMT} and \cite{Pau}), 
their use in probability theory seems to be quite recent. The CRT
is naturally
viewed as a random real tree (see Section 3 below), but this interpretation was not made
explicit in Aldous' work. The recent paper \cite{EPW} of Evans, Pitman and Winter starts 
with a general study of real trees from the point of view of measure theory, and 
establishes in particular that the space $\T$ of equivalent classes of (rooted) compact real trees,
endowed with the Gromov-Hausdorff metric, is a Polish space. This makes it very natural
to consider random variables or even random processes taking values in the space $\T$
(\cite{EPW} gives a particularly nice example of such a process, which combines root growth
and regrafting and converges in distribution to the CRT, see also \cite{EW} for further developments
along these lines). 

Our presentation owes a lot 
to the recent paper \cite{DuLG}, which uses the formalism of real trees to define the so-called L\'evy trees that
were implicit in \cite{LGLJ1} or \cite{DuLG0}, and to study various 
probabilistic and fractal properties of these objects.  In particular, the coding of a real tree by a
(deterministic) excursion is made precise in Section 2 below, which is 
taken from Section 2 of \cite{DuLG}. Aldous'
theorem relating the CRT to discrete Galton-Watson trees, and other analogous results involving  the more general
stable trees can be restated elegantly in the framework of real trees (see Sections 3 and 4 below). Stable trees
have a number of remarkable properties, some of which are briefly presented in Sections 4 and 5. 
In particular, they possess a nice self-similarity
property (Proposition 4.3) and they also verify a branching property
(Theorem \ref{existLT}) that is analogous
to the well-known branching property of Galton-Watson trees:
Conditionally given the tree
below level $a>0$, the subtrees originating from that level are distributed as the atoms of
a Poisson point measure whose intensity involves a local time measure supported on the 
vertices at distance $a$ from the root. These local times make it 
possible to define the uniform measure on the tree in an intrinsic way. 
Section 6 shows that the L\'evy snake 
construction of superprocesses, which generalizes the Brownian snake approach,
takes a neat form in the formalism of real trees (Theorem \ref{represuper}).
Finally, in Section 7, we give a very recent application of these concepts
to asymptotics for random quadrangulations, which involves a particular conditioning
of our random real trees.

\section{Real trees}

We start with a basic definition (see e.g. \cite{DMT}).

\begin{definition}
A metric space $(\t,d)$ is a real tree if the following two
properties hold for every $\sigma_1,\sigma_2\in \t$.

\begin{description}
\item{\rm(i)} There is a unique
isometric map
$f_{\sigma_1,\sigma_2}$ from $[0,d(\sigma_1,\sigma_2)]$ into $\t$ such
that $f_{\sigma_1,\sigma_2}(0)=\sigma_1$ and $f_{\sigma_1,\sigma_2}(
d(\sigma_1,\sigma_2))=\sigma_2$.
\item{\rm(ii)} If $q$ is a continuous injective map from $[0,1]$ into
$\t$, such that $q(0)=\sigma_1$ and $q(1)=\sigma_2$, we have
$$q([0,1])=f_{\sigma_1,\sigma_2}([0,d(\sigma_1,\sigma_2)]).$$
\end{description}

\noindent A rooted real tree is a real tree $(\t,d)$
with a distinguished vertex $\rho=\rho(\t)$ called the root.
\end{definition}

\noindent{\bf Remark} We use the terminology real tree rather than $\R$-tree as in \cite{DMT}
or \cite{Pau}.

\smallskip
In what follows, real trees will always be rooted and compact, even if this
is not mentioned explicitly. 

Let us consider a rooted real tree $(\t,d)$.
The range of the mapping $f_{\sigma_1,\sigma_2}$ in (i) is denoted by
$\llbracket \sigma_1,\sigma_2\rrbracket$ (this is the line segment between $\sigma_1$
and $\sigma_2$ in the tree). 
In particular, for every $\sigma\in \t$, $\llbracket \rho,\sigma\rrbracket$ is the path 
going from the root to $\sigma$, which we will interpret as the ancestral
line of vertex $\sigma$. More precisely we can define a partial order on the
tree by setting $\sigma\preccurlyeq \sigma'$
($\sigma$ is an ancestor of $\sigma'$) if and only if $\sigma\in\llbracket \rho,\sigma'
\rrbracket$.
If $\sigma,\sigma'\in\t$, there is a unique $\eta\in\t$ such that
$\llbracket \rho,\sigma
\rrbracket\cap \llbracket \rho,\sigma'
\rrbracket=\llbracket \rho,\eta
\rrbracket$. We write $\eta=\sigma\wedge \sigma'$ and call $\eta$ the most recent
common ancestor to $\sigma$ and $\sigma'$.

By definition, the multiplicity ${\bf k}(\sigma)$ of a vertex $\sigma\in\t$ is the number of
connected components of $\t\backslash \{\sigma\}$. Vertices
of $\t\backslash\{\rho\}$ which have multiplicity
$1$ are called leaves.

Let us now discuss convergence of real trees.
Two rooted real trees $\t_{(1)}$ and $\t_{(2)}$ are called equivalent if
there is a root-preserving isometry that maps 
$\t_{(1)}$ onto $\t_{(2)}$. We denote by ${\T}$ the set of all 
equivalence classes of rooted compact real trees. The set $\T$
can be equipped with the (pointed) Gromov-Hausdorff distance, which
is defined as follows. 

If $(E,\delta)$ is a metric space, we use the
notation $\delta_{Haus}(K,K')$ for the usual Hausdorff metric between
compact subsets of $E$.
Then, if $\t$ and $\t'$ are two rooted compact real trees
with respective roots $\rho$ and $\rho'$, we define the distance
$d_{GH}(\t,\t')$ as 
$$d_{GH}(\t,\t')=\inf\Big(\delta_{Haus}(\varphi(\t),\varphi'(\t'))\vee \delta(\varphi(\rho),
\varphi'(\rho'))\Big),$$ 
where the infimum is over all isometric embeddings $\varphi:\t\la E$ and
$\varphi':\t'\la E$ of $\t$ and $\t'$ into a common metric
space $(E,\delta)$. Obviously $d_{GH}(\t,\t')$ only depends on
the equivalence classes of $\t$ and $\t'$. Furthermore $d_{GH}$ defines 
a metric on $\T$ (cf \cite{Gro} and \cite{EPW}).

\begin{theorem}[\cite{EPW}]
\label{Polish}
The metric space $({\T},d_{GH})$
is complete and separable. 
\end{theorem}

Furthermore, the distance 
$d_{GH}$ can often be evaluated in the following way. First recall that if
$(E_1,d_1)$ and $(E_2,d_2)$ are two compact metric spaces, a
correspondence between $E_1$ and $E_2$ is a subset ${\cal R}$ of
$E_1\times E_2$ such that for every $x_1\in E_1$ there exists at least one
$x_2\in E_2$ such that $(x_1,x_2)\in{\cal R}$ and conversely
for every $y_2\in E_2$ there exists at least one
$y_1\in E_1$ such that $(y_1,y_2)\in{\cal R}$. The distorsion of
the correspondence ${\cal R}$ is defined by
$${\rm dis}({\cal R})=\sup\{|d_1(x_1,y_1)-d_2(x_2,y_2)|:
(x_1,x_2),(y_1,y_2)\in{\cal R}\}.$$
Then, if $\t$ and $\t'$ are two rooted $\R$-trees with respective roots
$\rho$ and $\rho'$, we have
\begin{equation}
\label{bounddist}
d_{GH}(\t,\t')={1\over 2}\ \inf_{{\cal R}\in{\cal
C}(\t,\t'),\,(\rho,\rho')\in{\cal R}}\,{\rm dis}({\cal R}),
\end{equation}
where ${\cal C}(\t,\t')$ denotes the set of all correspondences 
between $\t$ and $\t'$ (see Lemma 2.3 in \cite{EPW}).

\section{Coding compact real trees}

Our main goal in this section is to describe a method for constructing real trees, which
is particularly well-suited to our forthcoming applications to random trees.
We consider a	 (deterministic) continuous function
$g:[0,\infty)\longrightarrow[0,\infty)$ with compact support
and such that $g(0)=0$. To avoid
trivialities, we will also assume that $g$ is not identically zero.
For every $s,t\geq 0$, we set
$$m_g(s,t)=\inf_{r\in[s\wedge t,s\vee t]}g(r),$$
and
$$d_g(s,t)=g(s)+g(t)-2m_g(s,t).$$
Clearly $d_g(s,t)=d_g(t,s)$ and it is also easy to verify the triangle
inequality
$$d_g(s,u)\leq d_g(s,t)+d_g(t,u)$$
for every $s,t,u\geq 0$. We then introduce the equivalence relation
$s\sim t$ iff $d_g(s,t)=0$ (or equivalently iff $g(s)=g(t)=m_g(s,t)$). Let
$\t_g$ be the quotient space
$$\t_g=[0,\infty)/ \sim.$$
Obviously the function $d_g$ induces a distance on $\t_g$, and we keep the
notation $d_g$ for this distance. We denote by
$p_g:[0,\infty)\longrightarrow
\t_g$ the canonical projection. Clearly $p_g$ is continuous (when
$[0,\infty)$ is equipped with the Euclidean metric and $\t_g$ with the
metric $d_g$).

\begin{theorem}[\cite{DuLG}]
\label{tree-deterministic}
The metric space $(\t_g,d_g)$ is a real tree.
\end{theorem}

We will always view $(\t_g,d_g)$ as a rooted real tree with root $\rho=p_g(0)$.
If $\zeta>0$ is the supremum of the support of $g$, we have
$p_g(t)=\rho$ for every 
$t\geq \zeta$. In particular, $\t_g=p_g([0,\zeta])$ is compact.
We will call $\t_g$ {\it the real tree coded by $g$}.

A useful ingredient in the proof of Theorem \ref{tree-deterministic}
is the following root-change lemma, whose proof is an elementary
exercice.

\begin{lemma}
\label{root-change}
Let $s_0\in[0,\zeta)$. For any real $r\geq 0$, denote by
$\ov r$ the unique element of $[0,\zeta)$ such that
$r-\ov r$ is an integer multiple of $\zeta$. Set
$$g'(s)=g(s_0)+g(\ov{s_0+s})-2m_g(s_0,\ov{s_0+s}),$$
for every $s\in [0,\zeta]$, and $g'(s)=0$ for $s>\zeta$. Then,
the function $g'$ is continuous with compact support and satisfies
$g'(0)=0$, so that we can define
the metric space $(\t_{g'},d_{g'})$. Furthermore, for
every
$s,t\in[0,\zeta]$, we have
\begin{equation}
\label{iso-root}
d_{g'}(s,t)=d_g(\ov{s_0+s},\ov{s_0+t})
\end{equation}
and there exists a unique isometry $R$ from $\t_{g'}$ onto $\t_{g}$
such that, for every $s\in [0,\zeta]$,
\begin{equation}
\label{def-changeroot}
R(p_{g'}(s))=p_g(\ov{s_0+s}).
\end{equation}
\end{lemma}

Thanks to the lemma, the fact that $\t_g$  verifies property (i) in the 
definition of a real tree is obtained from the particular case when
$\sigma_1=\rho$ and $\sigma_2=\sigma=p_g(s)$ for some $s\in[0,\zeta]$.
In that case however, the isometric mapping $f_{\rho,\sigma}$ is easily constructed 
by setting
$$f_{\rho,\sigma}(t)=p_g(\sup\{r\leq s:g(r)=t\})\ ,\quad\hbox{for every }
0\leq t\leq g(s)=d_g(\rho,\sigma).$$
The remaining part of the argument is straightforward: See Section 2 in \cite{DuLG}
(or Lemma 3.1 in \cite{EW} for a different argument). 

\smallskip
\noindent{\bf Remark} In addition to the natural genealogical order, the tree
$\t_g$ is also equipped with the total order induced by the order
on the real line and the coding function $g$. We will not use this observation here
as we prefer to consider real trees as unordered trees (compare with the discrete
plane trees that are discussed in Section 3).

\smallskip
Can one compare the trees coded by two different functions $g$ and $g'$ ? The following
lemma gives a simple estimate.

\begin{lemma}
\label{dist-trees}
Let $g$ and $g'$ be two continuous functions with compact support
from $[0,\infty)$ into $[0,\infty)$, such that $g(0)=g'(0)=0$.
Then,
$$d_{GH}(\t_g,\t_{g'})\leq 2\|g-g'\|,$$
where $\|g-g'\|$ stands for the supremum norm of $g-g'$.
\end{lemma}

\noindent{\bf Proof.} We can construct a correspondence between 
$\t_g$ and $\t_{g'}$ by setting
$${\cal R}=\{(\sigma,\sigma'):\sigma=p_g(t)\hbox{ and }
\sigma'=p_{g'}(t)\hbox{ for some }t\geq 0\}.$$
In order to bound the distortion of $\cal R$, let
$(\sigma,\sigma')\in{\cal R}$ and $(\eta,\eta')\in{\cal R}$. By
our definition of ${\cal R}$ we can find $s,t\geq 0$ such that
$p_g(s)=\sigma$, $p_{g'}(s)=\sigma'$ and $p_g(t)=\eta$, $p_{g'}(t)=
\eta'$. Now recall that
$d_g(\sigma,\eta)=g(s)+g(t)-2m_g(s,t)$ and
$d_{g'}(\sigma',\eta')=g'(s)+g'(t)-2m_{g'}(s,t)$, 
so that
$$|d_g(\sigma,\eta)-d_{g'}(\sigma',\eta')|\leq 4\|g-g'\|.$$
Hence we have ${\rm dis}({\cal R})\leq 4\|g-g'\|$ and the desired result
follows from (\ref{bounddist}). \cq

The mapping $g\to \t_g$ is thus continuous if the set of  
functions satisfying our assumptions is equipped with the
supremum norm. In particular, this mapping is measurable.

\section{Galton-Watson trees and the CRT}

Denote by $({\bf e}_t)_{0\leq t\leq 1}$ a normalized Brownian excursion. Informally, $({\bf e}_t)_{0\leq
t\leq 1}$ is just a linear Brownian path started from the origin and conditioned to stay positive
over the time interval $(0,1)$, and to come back to $0$ at time $1$. See e.g. Sections 2.9 
and 2.10 of It\^o and McKean \cite{IM} for a discussion of the normalized excursion. We
extend 
the definition of ${\bf e}_t$ by setting ${\bf e}_t=0$ if $t>1$. Then the (random) function
${\bf e}$ satisfied the assumptions of the previous section and we can thus consider the
real tree $\t_{\bf e}$, which is a random variable with values in $\T$.

\begin{definition}
\label{CRTdef}
The random real tree $\t_{\bf e}$ is called the Continuum Random Tree (CRT).
\end{definition}

The CRT was initially defined by Aldous \cite{Al1} with a different formalism,
but the preceding definition corresponds to Corollary 22 in \cite{Al3}, up to an unimportant
scaling factor $2$.

One major motivation for studying the CRT is the fact that it occurs as the
scaling limit of 
critical Galton-Watson trees conditioned to have a large (fixed) number of vertices.
In order to state this result properly, we need to introduce some formalism for discrete trees,
which we borrow from Neveu \cite{Neveu}. Let
$${\cal U}=\bigcup_{n=0}^\infty \N^n  $$
where $\N=\{1,2,\ldots\}$ and by convention $\N^0=\{\varnothing\}$. If
$u=(u_1,\ldots u_m)$ and 
$v=(v_1,\ldots, v_n)$ belong to $\cal U$, we write $uv=(u_1,\ldots u_m,v_1,\ldots ,v_n)$
for the concatenation of $u$ and $v$. In particular $u\varnothing=\varnothing u=u$.

\begin{definition}
\label{plane-tree}
A plane tree $\theta$ is a finite subset of
$\cal U$ such that:
\begin{description}
\item{(i)} $\varnothing\in \theta$.

\item{(ii)} If $v\in \theta$ and $v=uj$ for some $u\in {\cal U}$ and
$j\in\N$, then $u\in\theta$.

\item{(iii)} For every $u\in\theta$, there exists a number $k_u(\theta)\geq 0$
such that $uj\in\theta$ if and only if $1\leq j\leq k_u(\theta)$.
\end{description}
\end{definition}

\noi We denote by ${\bf T}$ the set of all plane trees. In what follows, we see each vertex of the
tree $\theta$ as an individual of a population  whose $\theta$ is the family tree. 

If $\theta$ is a tree and $u\in \theta$, we define the shift of $\theta$ at $u$
by $\tau_u\theta=\{v\in {\cal U}:uv\in\theta\}$.
Note that $\tau_u\theta\in {\bf T}$. 
We also denote by $h(\theta)$ the height of $\theta$, that is the maximal generation of a vertex in $\theta$, and 
by $\#\theta$ the number of vertices of $\theta$.

\smallskip
For our purposes it will be convenient to view $\theta$ as a real tree: To this end,
embed $\theta$ in the plane, in such a way that each edge corresponds to
a line segment of length one, in the way suggested by the left part of Fig. 1. Denote by $\t^\theta$ the
union of all these line segments  and equip $\t^\theta$ with the obvious metric such that the distance
between 
$\sigma$ and
$\sigma'$ is the length of the shortest path from $\sigma$ to $\sigma'$ in $\t^\theta$. This construction
leads to a (compact rooted)  real tree whose
equivalence class does not depend on the particular embedding.

\bigskip
\begin{center}
\unitlength=1pt
\begin{picture}(320,130)

\thicklines \put(50,0){\line(-1,2){20}}
\thicklines \put(50,0){\line(1,2){20}}
\thicklines \put(30,40){\line(0,1){40}}
\thicklines \put(30,40){\line(-1,2){20}}
\thicklines \put(30,40){\line(1,2){20}}
\thicklines \put(30,80){\line(-1,2){20}}
\thicklines \put(30,80){\line(1,2){20}}

\thinlines \put(37,15){\vector(-1,2){6}}
\thinlines \put(17,55){\vector(-1,2){6}}
\thinlines \put(20,67){\vector(1,-2){6}}
\thinlines \put(27,57){\vector(0,1){12}}
\thinlines \put(17,95){\vector(-1,2){6}}
\thinlines \put(22,107){\vector(1,-2){6}}
\thinlines \put(33,95){\vector(1,2){6}}
\thinlines \put(48,107){\vector(-1,-2){6}}
\thinlines \put(33,69){\vector(0,-1){12}}
\thinlines \put(34,55){\vector(1,2){6}}
\thinlines \put(48,67){\vector(-1,-2){6}}
\thinlines \put(40,27){\vector(1,-2){6}}
\thinlines \put(54,15){\vector(1,2){6}}
\thinlines \put(68,27){\vector(-1,-2){6}}

\put(40,-5){$\varnothing$}
\put(20,35){$1$}
\put(60,35){$2$}
\put(12,75){\scriptsize{(1,1)}}
\put(30,75){\scriptsize{(1,2)}}
\put(50,75){\scriptsize{(1,3)}}
\put(12,115){\scriptsize{(1,2,1)}}
\put(50,115){\scriptsize{(1,2,2)}}

\thinlines \put(85,0){\vector(1,0){125}}
\thinlines \put(85,0){\vector(0,1){130}}
\thicklines \put(85,0){\line(1,5){8}}
\thicklines \put(93,40){\line(1,5){8}}
\thicklines \put(101,80){\line(1,-5){8}}
\thicklines \put(109,40){\line(1,5){8}}
\thicklines \put(117,80){\line(1,5){8}}
\thicklines \put(125,120){\line(1,-5){8}}
\thicklines \put(133,80){\line(1,5){8}}
\thicklines \put(141,120){\line(1,-5){8}}
\thicklines \put(149,80){\line(1,-5){8}}
\thicklines \put(157,40){\line(1,5){8}}
\thicklines \put(165,80){\line(1,-5){8}}
\thicklines \put(173,40){\line(1,-5){8}}
\thicklines \put(181,0){\line(1,5){8}}
\thicklines \put(189,40){\line(1,-5){8}}

\thinlines \put(93,0){\line(0,1){2}}
\thinlines \put(101,0){\line(0,1){2}}
\thinlines \put(109,0){\line(0,1){2}}
\thinlines \put(85,40){\line(1,0){2}}
\thinlines \put(85,80){\line(1,0){2}}
\thinlines \put(197,0){\line(0,1){2}}

\put(80,39){\scriptsize 1}
\put(80,79){\scriptsize 2}
\put(92,-6){\scriptsize 1}
\put(100,-6){\scriptsize 2}
\put(108,-6){\scriptsize 3}
\put(190,-6){$\scriptstyle 2(\#\theta -1)$}
\put (88,120){$C^\theta(t)$}
\put(205,4){$t$}

\thinlines\put(230,0){\vector(1,0){80}}
\thinlines\put(230,0){\vector(0,1){130}}
\thicklines \put(230,0){\line(1,4){10}}
\put(237.5,38){$\bullet$}
\thicklines \put(240,40){\line(1,4){10}}
\put(247.5,78){$\bullet$}
\thicklines\put(250,80){\line(1,0){10}}
\put(257.5,78){$\bullet$}
\thicklines\put(260,80){\line(1,4){10}}
\put(267.5,118){$\bullet$}
\thicklines\put(270,120){\line(1,0){10}}
\put(277.5,118){$\bullet$}
\thicklines\put(280,120){\line(1,-4){10}}
\put(287.5,78){$\bullet$}
\thicklines\put(290,80){\line(1,-4){10}}
\put(297.5,38){$\bullet$}

\thinlines \put(240,0){\line(0,1){2}}
\thinlines \put(250,0){\line(0,1){2}}
\thinlines \put(260,0){\line(0,1){2}}
\thinlines \put(300,0){\line(0,1){2}}
\thinlines \put(230,40){\line(1,0){2}}
\thinlines \put(230,80){\line(1,0){2}}

\put(239,-6){\scriptsize 1}
\put(249,-6){\scriptsize 2}
\put(259,-6){\scriptsize 3}
\put(290,-6){$\scriptstyle\#\theta-1$}
\put(305,4){$n$}
\put(233,120){$H^\theta_n$}

\put(225,39){\scriptsize 1}
\put(225,79){\scriptsize 2}

\put(30,-20){Tree $\theta$}
\put(100,-20){Contour function $C^\theta(t)$}
\put(230,-20){Height function $H^\theta_n$}

\end{picture}

\bigskip\medskip
\bigskip
Figure 1
\end{center}

\medskip

The discrete tree $\theta$, or equivalently the tree $\t^\theta$, can be coded by two
simple discrete functions, namely the {\it contour function} and the {\it height function}.
To define the contour function, consider a particle that starts from the root of $\t^\theta$
and visits continuously  all edges at speed one, going backwards as less as possible and respecting 
the lexicographical order of vertices. Then let $C^\theta(t)$ denote the distance to the root
of the position of the particle at time $t$ (for $t\geq 2(\#\theta-1)$, we take $C^\theta(t)=0$
by convention).
Fig.1 explains
the definition of the contour function better than a formal definition. Note that in the notation of
Section 2, we have $\t^\theta=\t_{C^\theta}$, meaning that $\t^\theta$ coincides with the tree
coded by the function $C^\theta$.

The height function $H^\theta=(H^\theta_n,0\leq n<\#\theta)$ is a discrete function defined as
follows. Write $u_0=\varnothing,u_1=1,u_2,\ldots,u_{\#\theta-1}$ for the elements of $\theta$
listed in lexicographical order. Then $H^\theta_n$ is the generation of $u_n$ (cf Fig.1).

Now let us turn to Galton-Watson trees.
Let $\mu$ be a critical offspring distribution.  This means that
$\mu$ is a probability measure on $\Z_+$ such that
$\sum_{k=0}^\infty k\mu(k)=1$.
We exclude the trivial case where $\mu(1)=1$.
Then,
there is a unique probability distribution $\Pi_\mu$ on ${\bf T}$ such that

\begin{description}
\item{(i)} $\Pi_\mu(k_\varnothing=j)=\mu(j)$,\quad
$j\in \Z_+$.

\item{(ii)} For every $j\geq 1$ with $\mu(j)>0$, the shifted trees
$\tau_1\theta,\ldots,\tau_j\theta$
are independent under the conditional probability
$\Pi_\mu(\cdot\mid k_\varnothing =j)$
and their conditional distribution is $\Pi_\mu$.
\end{description}

A random tree with distribution $\Pi_\mu$ is called a Galton-Watson tree
with offspring distribution $\mu$, or in short a $\mu$-Galton-Watson tree. Obviously
it describes the genealogy of the Galton-Watson process with offspring distribution $\mu$
started initially with one individual.

We are now able to state our invariance principle for Galton-Watson trees.
If $\t$ is a real tree with metric $d$ and $r>0$, we write $r\t$ for the  ``same'' tree
with metric $rd$. We say that $\mu$ is aperiodic if it is not supported on
a proper subgroup of $\Z$.

\begin{theorem}
\label{Aldous}
Assume that the offspring distribution $\mu$ is critical with finite 
variance $\sigma^2>0$, and is aperiodic. Then the
distribution of the rescaled tree ${\sigma\over 2\sqrt{n}}\t^\theta$ under the
probability measure $\Pi_\mu(\cdot \mid \#\theta=n)$ converges as $n\to \infty$
to the law of the CRT.
\end{theorem}

This is an immediate consequence of Theorem 23 in \cite{Al3}. In fact the latter result
states that the rescaled contour function $({\sigma\over 2\sqrt{n}}\,C^\theta(2nt),0\leq t\leq 1)$ under
$\Pi_\mu(\cdot \mid \#\theta=n)$ converges in distribution to the normalized
Brownian excursion $\bf e$. Note that the rescaled tree ${\sigma\over 2\sqrt{n}}\t^\theta$ is
coded by the function ${\sigma\over 2\sqrt{n}}\,C^\theta$, in the sense of Section 2 above. It then suffices to
apply Lemma \ref{dist-trees}.

A simple approach to the convergence of the rescaled contour functions towards the normalized
Brownian excursion was provided in \cite{MM0}. A key ingredient is the fact that the height function under $\Pi_\mu$
can be written as a simple functional of a random walk. Let 
$S$ be a random walk with jump distribution $\nu(i)=\mu(i+1)$ ($i=-1,0,1,2,\ldots$)
started from the origin, and $T=\inf\{n\geq 1:S_n=-1\}$. Then, under $\Pi_\mu$,
$$(H^\theta_n,0\leq n<\#\theta)\build{=}_{}^{\rm(d)}(K_n,0\leq n<T)$$
where
\begin{equation}
\label{keyform}
K_n={\rm Card}\Big\{0\leq j<n:S_j=\inf_{j\leq k\leq n} S_k\Big\}.
\end{equation}
Formula (\ref{keyform}) was noticed in \cite{LGLJ1}, where it motivated the coding of
generalizations of the CRT (see Section 4 below).  Under the additional assumption that
$\mu$ has exponential moments, formula (\ref{keyform}) was used in \cite{MM0} to show that
the height function (and then the contour function) under $\Pi_\mu(\cdot \mid \#\theta=n)$
is close to an excursion of the random walk $S$ of length $n$. Invariance principles
relating random walk excursions to Brownian excursions then lead to the desired result.

See also Bennies and Kersting \cite{BK} for a nice elementary presentation of the
relations between Galton-Watson trees and random walks, with an application to
a (weak) version of Theorem \ref{Aldous}, and Chapters 5 and 6 of Pitman \cite{Pit0}.

\smallskip
\noindent{\bf Remarks} (a) Theorem \ref{Aldous} has various applications to the
asymptotic behavior of functionals of Galton-Watson trees. For instance, by considering
the height of trees, we easily get that for every $x\geq 0$,
\begin{equation}
\label{height-limit}
\lim_{n\to\infty}
\Pi_\mu\Big(h(\theta)\geq {2x\sqrt{n}\over \sigma}\mid \#(\theta)=n\Big)=P(M({\bf e})\geq x),
\end{equation}
where
$$M({\bf e})=\sup_{s\geq 0} {\bf e}(s).$$
The probability in the right-hand side of (\ref{height-limit}) is known in the form of a series (Kennedy
\cite{Ke}, see also Chung \cite{Chung} for related results, and Section 3.1 in
\cite{Al2}).

\smallskip
\noindent(b) Special choices of the offspring distribution $\mu$ lead to limit theorems
for ``combinatorial trees''. For instance, if we let $\mu$ be the geometric distribution
$\mu(k)=2^{-k-1}$, which satisfies all our assumptions with $\sigma^2=2$, then 
$\Pi_\mu(\cdot\mid \#\theta=n)$ is easily identified as the uniform distribution on the
set of all plane trees with $n$ vertices. In particular, (\ref{height-limit}) gives the
asymptotic proportion of those plane trees with $n$ vertices which have a height greater
than $x\sqrt{2n}$. Similar observations apply to other classes of discrete trees, e.g.
to binary trees (take $\mu(0)=\mu(2)=1/2$, and note that we need a slight extension
of Theorem \ref{Aldous}, since $\mu$ is not aperiodic), or to Cayley trees (corresponding to
the case when $\mu$ is Poisson with parameter $1$). Such asymptotics had in fact been established
before Theorem \ref{Aldous} was proved, by the method of generating functions: See \cite{BKR} and especially
Flajolet and Odlyzko \cite{FO}.

\smallskip
\noindent(c) The convergence in Theorem \ref{Aldous} is not strong enough to 
allow one to deal with all interesting functionals of the tree. Still, Theorem 
\ref{Aldous} can be used to guess the kind of limit theorem one should expect.
A typical example is the height profile of the tree, that is for every level $k$
the number of vertices at generation $k$. Theorem \ref{Aldous} strongly suggests that
the suitably rescaled height profile of a Galton-Watson tree conditioned to have
$n$ vertices should converge in distribution towards the local time process of
a normalized Brownian excursion (see also the discussion of local times of 
stable trees in Section 5 below). This was indeed proved by Drmota and Gittenberger
\cite{DG} (see also Pitman \cite{Pit} and Aldous \cite{Al5}).

\smallskip
The effect of conditioning on the event $\{\#\theta=n\}$ is to force the tree $\theta$
to be large. One can imagine various other conditionings that have the same effect, and will
give rise to different limiting real trees. Typically, these limits will be described 
in terms of the It\^o excursion measure. Recall that the It\^o excursion measure $n(de)$ is the 
$\sigma$-finite measure on the space $C(\R_+,\R_+)$ of continuous functions from $\R_+$
and $\R_+$, which can be obtained as
\begin{equation}
\label{Ito-def}
n(de)=\lim_{\varepsilon\to 0} {1\over 2\varepsilon}\,P_\varepsilon(de)
\end{equation}
where $P_\varepsilon(de)$ stand for the distribution of a linear Brownian motion started
from $\varepsilon$ and stopped at the first time when it hits $0$, and 
we omit the precise meaning of the convergence (\ref{Ito-def}) (see Chapter XII of \cite{RY}
for a thorough discussion of It\^o's excursion measure). We write 
$\zeta(e)=\inf\{s>0:e(s)=0\}$ for the duration of $e$. 
The connection with the normalized Brownian excursion is made by noting that
the distribution of ${\bf e}$ is just the conditioned measure $n(de\mid \zeta(e)=1)$.

We write $\Theta_2(d\t)$ for the distribution of the tree $\t_e$ under $n(de)$ (the reason for the subscript $2$
will become clear later). 
Self-similarity properties of our random real trees become apparent
on the measure $\Theta_2$.

\begin{proposition}
\label{self-simi}
For every $r>0$, the distribution of $r\t$ under $\Theta_2$ is $r\Theta_2$. 
\end{proposition}

This result readily follows from the analogous statement for the It\^o measure.

To illustrate the usefulness of introducing the measure $\Theta$, we state a variant of Theorem
\ref{Aldous}. This variant and other results of the same type can be found in 
Chapter 2 of \cite{DuLG0} (see also Geiger and Kersting \cite{Gei2} for related results). The height
$h(\t)$ of a real tree is obviously defined as the maximal distance to the root.

\begin{theorem}
\label{Aldous-bis}
Under the same assumptions as in Theorem \ref{Aldous}, the distribution of the
tree $n^{-1/2}\t^\theta$ under $\Pi_\mu(d\theta \mid h(\theta)>n^{1/2})$ converges as
$n\to\infty$ to $\Theta_2(d\t \mid h(\t)>1)$. 
\end{theorem}

In contrast with Theorem \ref{Aldous}, this statement does not involve the constant $\sigma$.

\section{Stable trees}

One may ask what happens in Theorems \ref{Aldous} and \ref{Aldous-bis} in the case when
$\mu$ has infinite variance. Assuming that $\mu$ is in the domain of attraction
of a stable distribution with index $\alpha\in(1,2)$, we still get a limiting
random real tree, which is called the stable tree with index $\alpha$.

Before stating this result, we need to spend some time defining the limiting
object. We fix $\alpha\in(1,2)$ and consider a stable L\'evy process 
$X=(X_t)_{t\geq 0}$ with index $\alpha$ and no negative jumps, started 
from the origin. The absence of negative jumps
implies that $E[\exp(-\lambda X_t)]<\infty$ for every $\lambda\geq 0$, and we may normalize $X$
so that
$$E[\exp(-\lambda X_t)]=\exp(\lambda^\alpha t).$$

Set $I_t=\inf_{r\leq t} X_r$, which is a continuous process. Then it is well known that
the process $X-I$ is Markovian. Furthermore, the point $0$ is regular for $X-I$,
and the process $-I$ provides a local time at the origin for this Markov process. The associated
excursion measure is denoted by ${\bf N}$. It will play the same role as the
It\^o excursion measure (which one recovers in the case $\alpha=2$, up to inimportant scaling constants) in the 
previous paragraph. We again use the notation $\zeta$ for the duration of the excursion
under ${\bf N}$.

\begin{proposition}
\label{height-process}
There exists a continuous process $(H_s)_{s\geq 0}$, called the height process, such that
$$
H_s=\lim_{\varepsilon\to 0} {1\over \varepsilon}\int_0^s dr 
\;{\bf 1}\Big(X_r\leq \inf_{r\leq u\leq s} X_u+\varepsilon\Big)
\quad,\quad\hbox{for every }s\geq 0,\quad{\bf N}\ \hbox{a.e.}
$$
Moreover $H_0=0$ and $H_s=0$ for every $s\geq \zeta$, $\bf N$ a.e.
\end{proposition}

The formula for $H_s$ is a continuous analogue of the formula (\ref{keyform})
for the discrete height function of a Galton-Watson tree.

\begin{definition}
\label{stable-tree}
The measure $\Theta_\alpha$ defined as the law of the tree $\t_H$ under $\bf N$, is called the
distribution of the stable tree with index $\alpha$. We also 
denote by $\Theta^{(1)}_\alpha$ the distribution of $\t_H$
under the probability measure ${\bf N}(\cdot\mid \zeta=1)$.
\end{definition}

We can now state an analogue of Theorem \ref{Aldous} in the stable case.

\begin{theorem}[\cite{Duq}]
\label{stableGW}
Suppose that $\mu$ is critical and aperiodic, and that it is in the domain of attraction of
the stable law with index $\alpha\in(1,2)$, meaning that there exists a
sequence $a_n\uparrow\infty$ such that, if $Y_1,Y_2,\ldots$ are i.i.d. with distribution $\mu$,
$${1\over a_n}(Y_1+\cdots+Y_n-n)\build{\la}_{n\to\infty}^{\rm(d)} X_1.$$
Then, the law of the tree $n^{-1}a_n\,\t^\theta$ under $\Pi_\mu(\cdot\mid\#\theta=n)$
converges as $n\to \infty$ to $\Theta^{(1)}_\alpha$.
\end{theorem}

Other limit theorems relating discrete Galton-Watson trees to the measures $\Theta_\alpha$, in the spirit
of Theorem \ref{Aldous-bis}, can be found in Chapter 2 of \cite{DuLG0}
and in \cite{DuLG} (Theorem 4.1). These results apply more generally
to the L\'evy trees of \cite{DuLG}, which can be viewed as possible limits of sequences of critical Galton-Watson
trees  for which the offspring distribution depends on the tree taken in the sequence (in a way very similar
to the classical approximations of L\'evy processes by random walks). The recent article \cite{GK} gives a
limit theorem for Galton-Watson trees with possibly infinite variance, which is related to
the results of \cite{Duq} and \cite{DuLG0}.

The self-similarity property of Proposition \ref{self-simi} extends to stable trees:

\begin{proposition}
\label{self-simi-stable}
For every $r>0$, the distribution of $r\t$ under $\Theta_\alpha$ is $r^{{1\over \alpha-1}}\Theta_\alpha$. 
\end{proposition}

\section{Probabilistic and fractal properties of stable trees}

In this section we give some important properties of our stable trees. The following results, which are
taken from \cite{DuLG},
hold for any $\alpha\in(1,2]$.

Consider first a fixed real tree $(\t,d)$ in $\T$ with root $\rho(\t)$.
The level set of $\t$ at level $a>0$ is
$$\t(a)=\{\sigma\in\t:d(\rho(\t),\sigma)=a\}.$$
The truncation of the tree $\t$ at level $a$ is the new tree
$${\rm tr}_a(\t)=\{\sigma\in \t:d(\rho(\t),\sigma)\leq a\},$$
which is obviously equipped with the restriction of the distance $d$. It is easy to
verify that the mapping $\t\to {\rm tr}_a(\t)$
from $\T$ into itself is measurable.

Let us fix $a>0$ and denote by $\t^{(i),\circ}$, $i\in \ii$ the connected 
components of the open set
$$\t({(a,\infty)})=\{\sigma\in \t:d(\rho(\t),\sigma) >a\}.$$
Notice that the index set $\ii$ may be empty (if $h(\t)\leq a$), finite or
countable. Let $i\in \ii$. Then the ancestor of $\sigma$ at level $a$ must be the
same for every $\sigma\in \t^{(i),\circ}$. We denote by $\sigma_i$ this 
common ancestor and set $\t^{(i)}=\t^{(i),\circ}\cup\{\sigma_i\}$. Then $\t^{(i)}$ is a compact rooted
$\R$-tree with root
$\sigma_i$. The trees $\t^{(i)}$, $i\in I$ are called the subtrees
of $\t$ originating from level $a$, and we set
$$\n_a^\t:=\sum_{i\in \ii} \delta_{(\sigma_i,\t^{(i)})},$$
which is a point measure on $\t(a)\times \T$. 

\begin{theorem}
\label{existLT}
For every $a>0$ and for $\Theta_\alpha$ a.e. $\t\in\T$
we can define a finite measure $\ell^a$ on $\t$, in such a way that the
following properties hold:
\begin{description}
\item{\rm (i)} For every $a>0$, $\ell^a$ is 
supported on $\t(a)$, and $\{\ell^a\ne 0\}=\{h(\t)> a\}$, $\Theta_\alpha(d\t)$ a.e.
\item{\rm (ii)} For every $a>0$, we have $\Theta_\alpha(d\t)$ a.e. for every bounded
continuous function $\varphi$ on $\t$,
\begin{eqnarray}
\label{appLTlevel}
\langle\ell^a,\varphi\rangle&=&\lim_{\varepsilon\da 0}
C_\alpha\,\varepsilon^{{1\over \alpha-1}}\int \n^\t_a(d\sigma d\t')\,\varphi(\sigma)\,\un_{\{h(\t')\geq
\varepsilon\}}
\nonumber\\
&=&\lim_{\varepsilon\da 0}
C_\alpha\,\varepsilon^{{1\over \alpha-1}}\int \n^\t_{a-\varepsilon}(d\sigma
d\t')\,\varphi(\sigma)\,\un_{\{h(\t')\geq
\varepsilon\}}
\end{eqnarray}
\end{description}
where $C_\alpha=(\alpha-1)^{1/(\alpha-1)}$ if $1<\alpha<2$ and $C_2=2$.
Furthermore, for every $a>0$, the conditional distribution of the point measure
$\n^\t_a(d\sigma d\t')$, under the probability measure $\Theta_\alpha(d\t\mid h(\t)>a)$ and 
given ${\rm tr}_a(\t)$,
is that of a Poisson point measure
on $\t(a)\times \T$ with intensity $\ell^a(d\sigma)\Theta_\alpha(d\t')$.
\end{theorem}

The last property is the most important one. It may be called the {\it branching property}
of the stable tree as it is exactly analogous to the classical branching property for
Galton-Watson trees (cf Property (ii) in the definition of Galton-Watson trees in
Section 3). The random measure $\ell^a$ will be called the {\it local time}
of $\t$ at level $a$. Note that the normalization of local times, that is the particular choice of the constant
$C_\alpha$ in (\ref{appLTlevel}), is made so that the 
branching property holds in the form given in the theorem.

\smallskip
\noindent{\bf Remarks} (a) The branching property holds for the more general L\'evy trees which are
considered in \cite{DuLG}. Other classes of random real trees, which are related to various problems
of probability theory or combinatorics, are studied in \cite{AMP} and \cite{HaM}. A form of the
branching property was used by Miermont \cite{Mie} in his study of
fragmentations of the stable tree.

\smallskip
\noindent (b) It is possible to choose a modification of the collection $(\ell^a,a>0)$ which has good
continuity properties. Precisely, this process has a c\` adl\` ag modification, which is even
continuous when $\alpha=2$. When $\alpha<2$, the discontinuity points of the mapping
$a\to\ell^a$ correspond to levels of points of infinite multiplicity of the tree (cf 
Proposition \ref{multi} below). 

\smallskip
Local times can be used to give an intrinsic definition of the uniform measure
on the tree: We set
$${\bf m}=\int_0^\infty da\,\ell^a.$$
Note that if the tree is constructed as $\t=\t_H$ (resp. $\t=\t_{\bf e}$ when $\alpha=2$),
the measure ${\bf m}$ corresponds to the image measure of Lebesgue measure 
on $[0,\zeta]$ under the coding function $s\to p_H(s)$ (resp. $s\to p_{\bf e}(s)$).
Thanks to this observation, the probability measure $\Theta_\alpha^{(1)}$ can be 
identified with $\Theta_\alpha(\cdot\mid {\bf m}(\t)=1)$.

The next proposition gives precise information about the multiplicity
of vertices in our stable trees. Recall that ${\bf k}(\sigma)$ denotes the multiplicity of $\sigma$.

\begin{proposition}
\label{multi}
We have ${\bf m}(\{\sigma\in\t:{\bf k}(\sigma)>1\})=0$, $\Theta_\alpha$ a.e. (in other words 
almost all vertices are leaves). Moreover, we have  $\Theta_\alpha$ a.e.:
\begin{description}
\item{(i)} If $\alpha=2$, ${\bf k}(\sigma)\in\{1,2,3\}$ for all $\sigma\in\t$,
and the set $\{\sigma\in\t:{\bf k}(\sigma)=3\}$ of binary branching points is a countable dense 
subset of $\t$.
\item{(ii)} If $1<\alpha<2$, ${\bf k}(\sigma)\in\{1,2,\infty\}$ for all $\sigma\in\t$,
and the set $\{\sigma\in\t:{\bf k}(\sigma)=\infty\}$ of infinite branching points is a countable dense 
subset of $\t$.
\end{description}
\end{proposition}

Stable trees also enjoy a nice invariance property under uniform re-rooting, which is related
to the deterministic re-rooting Lemma \ref{root-change}. We state this
property under the law $\Theta_\alpha^{(1)}$ of the normalized tree. If $\t$ is a tree and $\sigma\in\t$, we write
$\t^{[\sigma]}$ for the ``same'' tree $\t$ with root
$\sigma$. 

\begin{proposition}
\label{re-rooting}
The law of the tree $\t^{[\sigma]}$ under the measure $\Theta_\alpha^{(1)}(d\t)\,{{\bf m}(d\sigma)}$
coincides with $\Theta_\alpha^{(1)}(d\t)$. 
\end{proposition}

In the case $\alpha=2$, this invariance property was already noticed in Aldous \cite{Al2}. 
Still for $\alpha=2$, more precise invariance properties under re-rooting can be found in
\cite{MM} and \cite{LGW}.

We conclude this section with a discussion of the fractal dimension of stable trees.
If $E$ is a subset of $\R_+$, we use the notation
$$\t(E)=\{\sigma\in\t:d(\rho,\sigma)\in E\}$$
where $\rho=\rho(\t)$ is the root of $\t$.

\begin{theorem}
\label{dimension}
Let $E$ be a compact subset of $(0,\infty)$ and $A=\sup E$. Assume that the Hausdorff and
upper box dimensions of $E$ are equal and let $d(E)\in [0,1] $
be their common value. Then, $\Theta_\alpha$ a.e. on the event $\{h(\t)>A\}$,
the Hausdorff and packing dimensions of $\t(E)$ coincide and are equal to
$${\rm dim}(\t(E))=d(E)+ {1\over \alpha-1}.$$
In particular,
$${\rm dim}(\t)={\alpha\over \alpha-1}\quad,\quad\Theta_\alpha\hbox{ a.e.}$$
and, for every $a>0$,
$${\rm dim}(\t(a))={1\over \alpha-1}\quad,\quad\Theta_\alpha\hbox{ a.e. on }\{h(\t)>a\}.$$
\end{theorem}

\noindent{\bf Remark} The formula for ${\rm dim}(\t)$ has also been derived by Haas 
and Miermont \cite{HaM} independently of \cite{DuLG}. More precise information about
the Hausdorff measure of stable trees can be found in \cite{DuLG2}.

\section{Spatial trees}

We will now explain how the genealogical structure of our stable trees
can be combined with spatial displacements (given by independent Brownian motions in $\R^k$) to yield a construction
of superprocesses with a stable branching mechanism. To present this construction in a
way suitable for applications, it is convenient to introduce the notion of a spatial tree.

Informally, a ($k$-dimensional) spatial tree is a pair $(\t,V)$ where 
$\t\in \T$ and $V$ is a continuous mapping from 
$\t$ into $\R^k$. Since we defined $\T$ as a space of equivalence
classes of trees, we should be a little more precise at this point.
If $\t$ and $\t'$ are two (rooted compact) real trees and $V$ and $V'$
are $\R^k$-valued continuous mappings defined respectively
on $\t$ and $\t'$, we say that the pairs $(\t,V)$ and $(\t,V')$ are
equivalent if there exists a root-preserving isometry $\Phi$
from $\t$ onto $\t'$ such that $V'_{\Phi(\sigma)}=V_\sigma$
for every $\sigma\in\t$. A spatial tree is then defined as an equivalent
class for the preceding equivalence relation, and we denote by $\T_{sp}$
the space of all spatial trees. Needless to say we will often abuse notation
and identify a spatial tree with an element of the corresponding equivalent class.

We denote by $\T_{sp}$ the set of all spatial trees. Recall the notation
of Section 1. We define a distance on $\T_{sp}$ by setting
$$d_{sp}((\t,V),(\t',V'))={1\over 2}\;\inf_{{\cal R}\in{\cal C}(\t,\t'),(\rho,\rho')\in{\cal R}}
\Big({\rm dis}({\cal R})+\sup_{(\sigma,\sigma')\in{\cal R}}|V_\sigma-V'_{\sigma'}|\Big),$$ where $\rho$ and $\rho'$
obviously denote the respective roots  of $\t$ and $\t'$. It is easy to verify that $(\T_{sp},d_{sp})$ is a Polish
space.

Let us fix $x\in\R^k$. Also let $\t\in\T$ be a compact rooted real tree 
with root $\rho$ and metric $d$. We may consider the $\R^k$-valued Gaussian process
$(Y_\sigma,\sigma\in
\t)$ whose distribution is characterized by
\begin{eqnarray*}
&&E[Y_\sigma]=x\;,\\
&&{\rm cov}(Y_\sigma,Y_{\sigma'})=d(\rho,\sigma\wedge \sigma')\,{\rm Id}\;,
\end{eqnarray*}
where ${\rm Id}$ denotes the $k$-dimensional identity matrix. The formula for the covariance is easy to
understand if we recall that $\sigma\wedge \sigma'$ is the most recent common ancestor to $\sigma$
and $\sigma'$, and so the ancestors of $\sigma$ and $\sigma'$ are the same up to
level $d(\varnothing,\sigma\wedge \sigma')$. Note that
$${\rm cov}(Y_\sigma-Y_{\sigma'},Y_\sigma-Y_{\sigma'})=d(\sigma,\sigma')\,{\rm Id}.$$
Let $\n(\t,\delta)$ stand for the minimal number of balls with radius $\delta$
needed to cover $\t$. From Theorem 11.17 in \cite{LT}, we
know that under the condition
\begin{equation}
\label{metricent}
\int_0^1 (\log \n(\t,\varepsilon^2))^{1/2}\;d\varepsilon<\infty,
\end{equation}
the process $(Y_\sigma,\sigma\in\t)$ has a continuous modification. We keep the notation
$Y$ for this modification. Assuming that (\ref{metricent}) holds, we denote by
$Q^x_\t$ the law on $\T_{sp}$ of $(\t,(Y_\sigma,\sigma\in\t))$.

It is easy to verify that condition (\ref{metricent}) holds
$\Theta_\alpha(d\t)$ a.e., for every $\alpha\in(1,2]$. The definition of $Q^x_\t$ then makes sense 
$\Theta_\alpha(d\t)$ a.e., and we may set
$$\N^{\alpha}_x=\int \Theta_\alpha(d\t)\,Q^x_\t,$$
which defines a $\sigma$-finite measure on $\T_{sp}$.

\smallskip

We can now turn to connections with superprocesses. Under the measure $\N^\alpha_x$, we may
for every $a>0$ define a measure $\z_a=\z_a(\t,V)$ on $\R^k$ by setting
$$\langle \z_a,\varphi\rangle=\int \ell^a(d\sigma)\,\varphi(V_\sigma),$$
where the local time measure $\ell^a$ was introduced in Theorem \ref{existLT}.
The next proposition reformulates a special case of Theorem 4.2.1 in \cite{DuLG}.

\begin{theorem} 
\label{represuper}
Let $\gamma\in M_f(\R^k)$ and let
$$\sum_{i\in \ii} \delta_{(\t^i,V^i)}$$
be a Poisson point measure on $\T_{sp}$ with intensity $\int \gamma(dx)\,\N^\alpha_x$. Then the 
process $(Z_a,a\geq 0)$ defined by
\begin{eqnarray*}
&&Z_0=\gamma\;,\\
&&Z_a=\sum_{i\in \ii} \z_a(\t^i,V^i)\;,\quad a>0\;,
\end{eqnarray*}
is a super-Brownian motion with branching mechanism $\psi(u)=u^\alpha$ 
($\psi(u)=2u^2$ if $\alpha=2$) started at $\gamma$. 
\end{theorem}

In the formula for $Z_a$, only finitely many terms can be nonzero, simply because
finitely many trees in the collection $(\t^i,i\in \ii)$ are such that $h(\t^i)>a$.
From the continuity properties of local times, we see that the version of $Z$ defined in the
proposition is c\` adl\` ag on $(0,\infty)$ for the weak topology on
finite measures on $\R^k$. By the known regularity properties of superprocesses
(see e.g. the more general Theorem 2.1.3 in \cite{DP}), it
must indeed be c\` adl\` ag on $[0,\infty)$. 

Theorem \ref{represuper} is clearly related to the Brownian snake construction of
superprocesses (in the case $\alpha=2$) or more generally to the L\'evy snake of
\cite{LGLJ2} or \cite{DuLG0}. Consider the case $\alpha=2$ and assume that
the spatial displacements have constructed (with initial point $x$)
using the tree $\t_{e}$ (where $e$ is distributed according to
the It\^o measure). For every $s\geq 0$, we may consider the path 
consisting of the spatial positions $V_\sigma$ along the ancestor line 
of the vertex $p_{e}(s)$. This gives a random path $W_s$ with duration
$d_{e}(\rho,p_{e}(s))={e}(s)$. The process $(W_s)_{s\geq 0}$ is then 
the path-valued Markov process called
the Brownian snake, here constructed under its excursion measure from $x$.

In view of Theorem \ref{represuper}, the measures $\N^\alpha_x$ (or rather the distribution
under $\N^\alpha_x$ of the measure-valued process $(\z_a,a\geq 0)$) are called the excursion 
measures of the $\psi$-super-Brownian motion. In the quadratic branching case, these
measures play an important role in the study of connections between superprocesses
and partial differential equations: See in particular \cite{LG1}. In the case
of a general branching branching mechanism, excursion measures are constructed via
the L\'evy snake in Chapter 4 of \cite{DuLG}, and a different approach has been 
proposed recently by Dynkin and Kuznetsov \cite{DK}.

Connections with partial differential equations are helpful to
analyse the probabilistic properties of spatial trees. We content ourselves
with one statement concerning the range 
$${\cal R}=\{V_\sigma:\sigma\in \t\}.$$

\begin{theorem}
\label{hitting}
Let $K$ be a compact subset of $\R^k$. Then the function
$$u(x)=\N^\alpha_x({\cal R}\cap K\not =\emptyset)\ ,\ x\in\R^k\backslash K$$
is the maximal nonnegative solution of ${1\over 2}\Delta u=u^\alpha$
(${1\over 2}\Delta u=2u^2$ if $\alpha=2$)
in $\R^k\backslash K$.
\end{theorem}

This theorem is a reformulation in our formalism of results proved by Dynkin \cite{Dy0}
in the framework of the theory of superprocesses. 

\smallskip
\noindent{\bf Integrated super-Brownian excursion (ISE)}.
Let us return for a while to the CRT, and combine the branching structure 
of the CRT with $k$-dimensional Brownian motions started from $x=0$, in the
way explained in the previous section. Precisely this means that we are
considering the probability measure on $\T_{sp}$ defined by
$$\N^{2,(1)}_0=\int \Theta_2^{(1)}(d\t)\,Q^0_\t,$$
where $\Theta_2^{(1)}(d\t)$ is the law of the CRT, in agreement with our
previous notation $\Theta_\alpha^{(1)}$.
Recall the notation ${\bf m}$ for the uniform measure on $\t$
(this makes sense $\Theta_2^{(1)}(d\t)$ a.s.).

\begin{definition}
The random probability measure $U$ on $\R^k$ defined under $\N^{2,(1)}_0$ by 
$$\langle U,\varphi\rangle = \int {\bf m}(d\sigma)\,\varphi(V_\sigma)$$
is called $k$-dimensional ISE (for Integrated Super-Brownian Excursion).
\end{definition}

Note that the topological support of ISE is the range ${\cal R}$ 
of the spatial tree, and that ISE should be interpreted as the
uniform measure on this set.
The random measure
ISE was first discussed by Aldous \cite{Al4}. It occurs in various 
asymptotics for models of statistical mechanics:
See in particular \cite{Sl1} and \cite{Sl2}.

\section{Conditioned spatial trees and quadrangulations}

In this section, we study a conditioning problem for the spatial
trees of the previous section, and we then explain why this conditioning
problem is related to certain asymptotics for planar quadrangulations.
We consider {\bf one-dimensional} spatial displacements (${k=1}$) and we write $\N^{(1)}_0$ for $\N^{2,(1)}_0$ since
we will only consider the case $\alpha=2$. Recall the notation $(\t,(V_\sigma)_{\sigma\in\t})$
for the generic element of $\T_{sp}$, and $\r$ for the associated range.

Our goal is to define the probability measure $\N^{(1)}_0$ 
conditioned on the event $\{\r\subset [0,\infty)\}$. In other
words we consider a tree of Brownian paths started from the
origin, whose underlying genealogical structure is the CRT,
and we want to condition this tree to remain on the positive half-line. Obviously this
conditioning is very degenerate (it is already degenerate if one considers a 
single Brownian path started from the origin). 

\begin{theorem}[\cite{LGW}]
\label{main-cond}
We have 
$$\lim_{\varepsilon\da
0}\varepsilon^{-4}\,\N^{(1)}_0(\r\subset(-\varepsilon,\infty))={2\over 21}.$$
There
exists a probability measure on $\T_{sp}$, which is denoted by $\ov\N^{(1)}_0$,
such that
$$\lim_{\varepsilon\da 0} \N^{(1)}_0(\cdot\mid
\r\subset(-\varepsilon,\infty))=\ov\N^{(1)}_0,$$ in the sense of weak convergence in
the space of probability measures on
$\T_{sp}$.
\end{theorem}

The measure $\ov\N^{(1)}_0$ is the law of the conditioned spatial tree we were
aiming at. The next result will show that, rather remarkably, this conditioned tree can be obtained
from the unconditioned one simply by re-rooting it at the vertex 
corresponding to the
minimal spatial position.

Before stating this result we need some notation for re-rooted trees. If 
$(\t,V)$ is a spatial tree, and $\sigma_0\in\t$, the re-rooted spatial
tree $(\t^{[\sigma_0]},V^{[\sigma_0]})$ is defined by saying that 
$\t^{[\sigma_0]}$ is $\t$ re-rooted at $\sigma_0$, and $V^{[\sigma_0]}_\sigma
=V_{\sigma}-V_{\sigma_0}$ for every $\sigma\in\t$.

\begin{theorem}[\cite{LGW}]
\label{verwaat}
There is $\N^{(1)}_0$ a.s. a unique vertex $\sigma_*$ 
which minimizes $V_\sigma$ over $\sigma\in\t$.
The probability measure $\ov\N^{(1)}_0$ is the law under 
$\N^{(1)}_0$ of the re-rooted spatial tree $(\t^{[\sigma_*]},V^{[\sigma_*]})$.
\end{theorem}

This theorem is reminiscent of a famous
theorem of Verwaat \cite{Verwaat} connecting the Brownian bridge and the
normalized Brownian excursion. The proof of Theorem \ref{verwaat} makes a heavy use
of the invariance property under re-rooting (Proposition \ref{re-rooting}).

It is interesting to reinterpret the preceding theorem in terms of ISE.
If we want to define
one-dimensional ISE conditioned to put no mass on the negative half-line,
a natural way is to condition it to put no mass on $]-\infty,-\varepsilon[$ 
and then to let $\varepsilon$ go to $0$. As a consequence of 
the previous two theorems, this is equivalent to shifting
the unconditioned ISE to the right, so that the left-most point
of its support becomes the origin. 

Our motivation for studying conditioned spatial trees came from applications
to planar quadrangulations. In order to describe these 
applications, we first need to introduce well-labelled trees.
We call labelled tree any pair $(\theta,v)$ where $\theta\in {\T}$
is a plane tree (cf Section 3) and $v$ is a mapping from $\theta$
into the set $\Z$ of integers. This is the obvious discrete analogue of the spatial trees considered
above. We say that $(\theta,v)$ is a well-labelled tree if in addition $v(\varnothing)=1$,
$v(u)\geq 1$ for every $u\in\theta$ and $|v(u)-v(u')|\leq 1$ whenever 
$u$ and $u'$ are neighboring vertices (that is, $u$ is the father of $u'$ or
$u'$ is the father of $u$). Again, this is a discrete version of the conditioned
spatial trees discussed in the preceding theorems.

In a way similar to what we did in Section 3, there is an obvious way
of viewing a labelled tree as a spatial tree: If $(\theta,v)$ is a labelled tree,
let $\t^\theta$ be as in Section 3, and define $V^{(\theta,v)}$ by setting
$V^{(\theta,v)}(\sigma)=v(u)$ if $\sigma$ is the vertex of $\t^\theta$
corresponding to $u$, and then interpolating linearly between neighboring 
vertices to complete the definition of $V^{(\theta,v)}$.

Let ${\cal W}_n$ stand for the set of all well-labelled trees with 
$n$ edges (or equivalently $n+1$ vertices).

\begin{theorem}
\label{invarprin}
The law of the rescaled tree $((2n)^{-1/2}\t^\theta,(9/8)^{1/4}n^{-1/4}V^{(\theta,v)})$
under the uniform probability measure on ${\cal W}_n$ converges as
$n\to\infty$ towards the measure $\ov\N^{(1)}_0$.
\end{theorem}

This theorem is a consequence of more general statements obtained in \cite{LG5}
for conditioned Galton-Watson trees. Closely related results can be found 
in \cite{CS} and \cite{MM}. The factor $(9/8)^{1/4}$ is easy to understand if
we write $(9/8)^{1/4}=2^{-1/4}(2/3)^{-1/2}$ and note that $2/3$ is the variance of the
uniform distribution on $\{-1,0,1\}$ (while the factor $2^{-1/4}$
corresponds to the term $2^{-1/2}$ in $(2n)^{-1/2}\t^\theta$).

\smallskip
\noindent{\bf Remark} Without the positivity constraint (for instance, considering labelled trees
with the only properties that the label of the root is $0$ and the labels of two neighboring
vertices differ by at most $1$), the limiting distribution would be 
$\N^{(1)}_0$ instead of $\ov\N^{(1)}_0$. This unconditional analogue of Theorem \ref{invarprin} follows
rather easily from Theorem \ref{Aldous}: See Janson and Marckert \cite{JM} for much
more general statements of this type (Kesten \cite{Kes} discusses related results
for tree-indexed random walks under different conditionings of the underlying trees).

\smallskip

Let us now discuss quadrangulations.
A planar quadrangulation is a planar map where
each face, including the unbounded one, has degree $4$. We are interested in
rooted quadrangulations, meaning that we distinguish an edge on the border
of the infinite face, which is 
oriented counterclockwise and called the root edge. The origin of the root edge is called
the root vertex. Two rooted quadrangulations are considered identical if there 
is a homeomorphism of the plane that sends one map onto the other.
We refer to Chassaing and Schaeffer
\cite{CS} for more precise definitions. We denote by ${\cal Q}_n$
the set of all rooted quadrangulations with $n$ faces. 
A key result gives a bijection
between the sets ${\cal W}_n$ and ${\cal Q}_n$
(see Theorem 1 in \cite{CS}, and note that this bijection has been extended recently 
to more general planar maps by Bouttier, Di Francesco and Guitter \cite{BDG3}). Moreover the radius of the
quadrangulation, defined as the maximal graph distance between the root vertex
and another vertex, corresponds via this bijection to the maximal label
of the tree. We can then deduce the following corollary (\cite{CS}, Corollary 3).

\begin{corollary}
\label{radius}
Let $R_n$ denote the radius of a random quadrangulation chosen uniformly
among all rooted planar quadrangulations with $n$ faces. Then,
$$n^{-1/4}R_n\build{\la}_{n\to\infty}^{\rm(d)} \rho$$
where the limiting variable $\rho$ is distributed as 
$$\Big({8\over 9}\Big)^{1/4}\Big(\sup_{\sigma\in \t} V_\sigma -\inf_{\sigma\in \t} V_\sigma\Big)$$
under $\N^{(1)}_0$ (alternatively, $\rho$ is distributed as $(8/9)^{1/4}$ times the length of the support
of one-dimensional ISE).
\end{corollary}

The proof of the corollary is easy from the preceding observations. Properties 
of the bijection between the sets ${\cal W}_n$ and ${\cal Q}_n$ imply that
$R_n$ has the same distribution as
$\sup\{v(u):u\in \theta\}$ under the uniform probability measure on ${\cal W}_n$.
By construction,
$$n^{-1/4}\sup\{v(u):u\in \theta\}=n^{-1/4}\sup\{V^{(\theta,v)}(\sigma):\sigma\in\t^\theta\}.$$
By Theorem \ref{invarprin}, the law of the latter quantity under 
the uniform probability measure on ${\cal W}_n$ converges to the law of
$$\Big({8\over 9}\Big)^{1/4}\;\sup_{\sigma\in\t} V_\sigma$$
under $\ov\N^{(1)}_0$. By Theorem \ref{verwaat}, this distribution is the same 
as the limiting one in Corollary
\ref{radius}. Note some information about this limiting distribution can be
found in Delmas \cite{Delmas} and in the recent preprint \cite{Bo}.

\smallskip
Large random quadrangulations, or more general planar maps, are used in theoretical physics
as models of random surfaces (see in particular \cite{BDG}, \cite{BDG2} and \cite{Dur}). Assuming that
there is a limiting continuous object for uniform random quadrangulations with $n$ faces, the normalizing factor
$n^{-1/4}$ in Corollary \ref{radius} suggests that its fractal dimension should be $4$, a fact
that is widely believed in the physics literature. A recent paper of Marckert and
Mokkadem \cite{MM} uses ideas related to the previous discussion 
to construct the Brownian map, which is a candidate for the continuous
limit of random quadrangulations.

\medskip
\noindent{\bf Acknowledgement} The authors wishes to thank the referee for several remarks and
for providing additional references.

\end{document}